\documentclass[11pt,twoside]{article}

\title{Concentration in Lotka-Volterra parabolic or integral equations:\\ a general convergence result}
\author{Guy Barles\thanks{Laboratoire de Math\'ematiques et Physique Th\'eorique, CNRS UMR 6083, F\'ed\'eration Denis Poisson,
Universit\'e Fran\c{c}ois Rabelais, Parc de Grandmont,
37200 Tours, France.
Email: barles@lmpt.univ-tours.fr}
, Sepideh Mirrahimi\thanks{
Universit\'e Pierre et Marie Curie-Paris 6, UMR 7598 LJLL, BC187, 4, place Jussieu,  F-75252 Paris cedex 5. Email: mirrahimi@ann.jussieu.fr}
, Beno\^{\i}t Perthame \footnotemark[2] \thanks{Institut Universitaire de France. Email: benoit.perthame@upmc.fr}
}

\usepackage[psamsfonts]{amssymb}

\usepackage{latexsym,amsfonts,graphicx}

\usepackage{amscd,amsmath,makeidx}

\usepackage{amsthm}
\usepackage{latexsym}
\usepackage[ansinew]{inputenc}
\usepackage{exscale}
\usepackage{dsfont} %

\def\bR{{\mathbb R}}

\def\cC{{\cal C}}

\def\rC{\mathrm{C}}
\def\ignore#1{}
\newtheorem{satz}{Satz}[section]

\newtheorem{theorem}[satz]{Theorem}
\newtheorem{lemma}[satz]{Lemma}

\newcommand{\e}{\epsilon}
\newcommand{\f}{\frac}
\newcommand{\p}{\partial}

\newcommand{\R}{\mathbb{R}}

\begin{document}

\maketitle
\pagestyle{plain}
%\tableofcontents
\pagenumbering{arabic}

\begin{abstract}
We study two equations of Lotka-Volterra type that describe the Darwinian evolution of a population density. In the first model a Laplace term represents the mutations. In the second one we model the mutations by an integral kernel. In both cases, we use a nonlinear birth-death term that corresponds to the competition between the traits leading to selection. \\
In the limit of rare or small mutations, we prove that the solution converges to a sum of moving Dirac masses. This limit is described by a constrained Hamilton-Jacobi equation. This was already proved in \cite{b.a} for the case with a Laplace term. Here we generalize the assumptions on the initial data and prove the same result for the integro-differential equation.
\end{abstract}


\noindent{\bf AMS Class. No:} 35B25, 35K57, 47G20, 49L25, 92D15

\section{Introduction}
We continue the study, initiated in \cite{b.a}, of the asymptotic behavior of Lotka-Volterra parabolic equations. The model we use describes the dynamics of a population density. Individuals respond differently to the environment, i.e. they have different abilities to use the available resources. To take this fact into account, population models can be structured by a parameter, representing a physiological (phenotypical) trait inherited from the parent, and that we denote by $x\in \bR^d$. We denote by $n(t,x)$ the density of trait $x$. The mathematical modeling in accordance with Darwin's theory consists of two effects: natural selection and mutations between the traits (see \cite{b.v, b.y, b.x,b.bb} for literature in adaptive evolution). We represent the birth and death rates of the phenotypical traits bya net growth rate $R(x,I)$. The term $I(t)$ is an ecological parameter that corresponds to a measure of the total population, whatever the trait, and that represents in the simpler possible way the resources (more precisely the inverse of it). We use two different models for mutations. A first possibility is to represent them by a Laplacian and, in an extreme and irrealistic simplification, we take them independent of birth, so as to write

\begin{equation}
\label{eq.1}
\begin{cases}
\partial_t n_\epsilon-\epsilon \triangle n_\epsilon=\frac{n_\epsilon}{\epsilon}R(x,I_\epsilon(t)),\;\;\;\;\; x\in \bR^d,\,t\geq0,\\
n_\epsilon(t=0)=n_\epsilon^0\in L^1(\mathbb{R}^d),\;\;\;\;\; n_\epsilon^0\geq0,
\end{cases}
\end{equation}

\begin{equation}
\label{eq.2}
I_\epsilon(t)=\int_{\bR^d}\psi(x)\,n_\epsilon(t,x)dx.
\end{equation}
Here $\epsilon$ is a small term that we introduce to consider only rare mutations. It is also used to re-scale time to consider a much larger time than a generation scale.\\

A more natural way to model mutations is to use, instead of a Laplacian, an integral term that describes directly the mutation probability to generate a new-born of trait $x$ from a mother with trait $y$. This yields 

\begin{equation}
\label{eq.i1}
\begin{cases}
\partial_t n_\epsilon=\frac{n_\epsilon}{\epsilon}R(x,I_\epsilon(t))+\frac{1}{\epsilon}\int \frac{1}{\epsilon^d}K(\frac{y-x}{\epsilon})\, b(y,I_\epsilon) \,n_\epsilon(t,y)\,dy,\;\;\;\;\; x\in \bR^d,\,t\geq0,\\
n_\epsilon(t=0)=n_\epsilon^0\in L^1(\mathbb{R}^d),\;\;\;\;\; n_\epsilon^0\geq0,
\end{cases}
\end{equation}

\begin{equation}
\label{eq.i2}
I_\epsilon(t)=\int_{\bR^d} n_\epsilon(t,x)dx.
\end{equation}

Both types of models can be derived from individual based stochastic processes in the limit of large populations depending on the scales in mutations birth and death (see \cite{b.aa,b.z}).
\\

In this paper, we study the asymptotic behavior of equations (\ref{eq.1})-(\ref{eq.2}) and (\ref{eq.i1})-(\ref{eq.i2}) when $\epsilon$ vanishes. Our purpose is to show that under some assumptions on $R(x,I)$, $n_\epsilon(t,x)$ concentrates as a sum of Dirac masses that are traveling. In biological terms, at every moment one or several dominant traits coexist while other traits disappear. The dominant traits change in time due to the presence of mutations.\\

We use the same assumptions as \cite{b.a}. We assume that there exist two constants $\psi_m$, $\psi_M$ such that

\begin{equation}
\label{eq.3}
0<\psi_m<\psi<\psi_M<\infty, \hspace{20 pt} \psi\in W^{2, \infty}(\bR^d).
\end{equation}

We also assume that there are two constants $0<I_m<I_M<\infty$ such that
\begin{equation}
\label{eq.4}
\min_{x\in \bR^d} R(x,I_m)=0, \hspace{20 pt} \max_{x\in \bR^d} R(x,I_M)=0,
\end{equation}
and there exists constants $K_i>0$ such that, for any $x\in \bR^d$, $I\in \bR$,
\begin{equation}
\label{eq.5}
-K_1\leq \frac{\partial R}{\partial I}(x,I)\leq-K_1^{-1}<0,
\end{equation}

\begin{equation}
\label{eq.5b}
\sup_{\frac{I_m}{2}\leq I \leq 2I_M} \parallel R(\cdot,I)\parallel_{W^{2,\infty}(\bR^d)}<K_2.
\end{equation}

We also make the following assumptions on the initial data

\begin{equation}
\label{eq.6b}
\hspace{10 pt}I_m \leq \int_{\bR^d}\psi(x)n_\epsilon^0(x)\leq I_M,\;\;\; \text{and} \;\;\;\exists\, A,\,B>0\,,\; n_\epsilon^0\leq e^{\frac{-A|x|+B}{\epsilon}},
\end{equation}
and that there exist a point $x_0\in \R^d$ and positive constants $L_0$ and $M_0$ such that
\begin{equation}
\label{as:n0-below}
 e^{-\f {M_0}{ \e}} \leq n_\e^0(x),\quad \text{for all $|x-x_0|\leq L_0$},
\end{equation}
Note that assumption (\ref{as:n0-below}) means that initially  we have some kind of biodiversity since it can be interpreted as different traits being sufficiently represented in the population. \\

Here we take $\psi(x)\equiv 1$ for equations (\ref{eq.i1})-(\ref{eq.i2}) because replacing $n$ by $\psi n$ leaves the model unchanged.
For equation (\ref{eq.i1}) we assume additionally that the probability kernel $K(z)$ and the mutation birth rate $b(z)$ verify
\begin{equation}
\label{eq.5i}
0\leq K(z),\;\;\;\;\; \int K(z)\, dz =1, \;\;\;\;\; \int K(z)e^{|z|^2}\, dz<\infty,
\end{equation}
\begin{equation}
\label{eq.6i}
b_m\leq b(z,I)\leq b_M,\;\;\;\;|\nabla_x b(z,I)|< L_1\,b(z,I),\;\;\;\;|b(x,I_1)-b(x,I_2)|<L_2|I_1-I_2|,
\end{equation}
where $b_m$, $b_M$, $L_1$ and $L_2$ are positive constants. Finally for equation (\ref{eq.i1}) we replace (\ref{eq.4}) and (\ref{eq.5}) by
\begin{equation}
\label{eq.6ii}
\min_{x\in \bR^d} \big[R(x,I_m)+b(x,I_m)\big]=0, \hspace{20 pt} \max_{x\in \bR^d} \big[R(x,I_M)+b(x,I_M)\big]=0,
\end{equation}
\begin{equation}
\label{eq.6ij}
|R(x,I_1)-R(x,I_2)|<K_3|I_1-I_2|\;\;\;\;\text{and}\;\;\;\;-K_4\leq \frac{\partial (R+b)}{\partial I}(x,I)\leq-K_4^{-1}<0,
\end{equation}
where $K_3$ and $K_4$ are positive constants. \\

In both cases, in the limit we expect $n(t,x)=0$ or $R(x,I)=0$, where $n(t,x)$ is the weak limit of $n_\epsilon(t,x)$ as $\epsilon$ vanishes. If we suppose that the latter is possible at only isolated points, we expect $n$ to concentrate as Dirac masses. Following earlier works on the similar issue \cite{b.e, b.h,b.a,b.w}, in order to study $n$, we make a change of variable $n_\epsilon(t,x)=e^\frac{u_\epsilon(t,x)}{\epsilon}$. It is easier to study the asymptotic behavior of $u_\epsilon$ instead of $n_\epsilon$. In section \ref{asymp} we study the asymptotic behavior of $u_\epsilon$ while $\epsilon$ vanishes. We show that $u_\epsilon$, after extraction of a subsequence, converge to a function $u$ that satisfies a constrained Hamilton-Jacobi equation in the viscosity sense (see \cite{b.d,b.n,b.u, b.dd} for general introduction to the theory of viscosity solutions). Our main results are as follows.

\begin{theorem}
\label{th.3}
Assume (\ref{eq.3})-(\ref{as:n0-below}). Let $n_\epsilon$ be the solution of (\ref{eq.1})-(\ref{eq.2}), and $u_\epsilon=\epsilon \ln(n_\epsilon)$. Then, after extraction of a subsequence, $u_\epsilon$ converges locally uniformly to a function $u\in \rC((0,\infty)\times \bR^d)$, a viscosity solution to the following equation:

\begin{equation}
\label{as.2}
\begin{cases}
\partial_t u=|\nabla u|^2+R(x,I(t)),\\
\underset{x\in \bR^d}{\max} \;u(t,x)=0,\;\;\;\forall t>0,
\end{cases}
\end{equation}

\begin{equation}
\label{as.1.1}
I_\epsilon(t) \underset{\epsilon \rightarrow 0}\longrightarrow I(t)\;\;\;\; \text{a.e.},\;\;\;\; \int \psi(x) n(t,x)dx=I(t)\;\;\;\;\text{a.e.}.
\end{equation}

In particular, a.e. in $t$, $supp\; n(t,\cdot)\subset \{u(t,\cdot)=0\}$. Here the measure $n$ is the weak limit of $n_\epsilon$ as $\epsilon$ vanishes. If additionally $(u_\epsilon^0)_\epsilon:= \epsilon \ln(n_\epsilon^0)$ is a sequence of uniformly continuous functions which converges locally uniformly to $u^0$ then $u\in \rC([0, \infty)\times \bR^d)$ and $u(0,x)=u^0(x)$ in $\bR^d$.
\end{theorem}

\begin{theorem}
\label{th.i2}
Assume (\ref{eq.5b})-(\ref{eq.6ij}), and $(u_\epsilon^0)_\epsilon$ is a sequence of uniformly Lipschitz-continuous  functions which converges locally uniformly to $u^0$. Let $n_\epsilon$ be the solution of (\ref{eq.i1})-(\ref{eq.i2}) with $n_\epsilon^0=e^{\frac{u\epsilon^0}{\epsilon}}$, and $u_\epsilon=\epsilon \ln(n_\epsilon)$. Then, after extraction of a subsequence, $u_\epsilon$ converges locally uniformly to a function $u\in \rC([0, \infty)\times \bR^d)$, a viscosity solution to the following equation:

\begin{equation}
\label{as.i1}
\begin{cases}
\partial_t u=R(x,I(t))+b(x,I(t))\int K(z)e^{\nabla u\cdot z}dz,\\
\underset{x\in \bR^d}{\max} \;u(t,x)=0,\;\;\;\forall t>0,\\
u(0,x)=u^0(x),
\end{cases}
\end{equation}

\begin{equation}
\label{as.i1.1}
I_\epsilon(t) \underset{\epsilon \rightarrow 0}\longrightarrow I(t)\;\;\;\; \text{a.e.},\;\;\;\; \int n(t,x)dx=I(t)\;\;\;\;\text{a.e.}.
\end{equation}

In particular, a.e. in $t$, $supp\; n(t,\cdot)\subset \{u(t,\cdot)=0\}$. As above, the measure $n$ is the weak limit of $n_\epsilon$ as $\epsilon$ vanishes.

\end{theorem}

These theorems improve previous results proved in \cite{b.e,b.a,b.h,b.l} in various directions. For the case where mutations are described by a Laplace equation, i.e. (\ref{eq.1})-(\ref{eq.2}), Theorem \ref{th.3} generalizes the assumptions on the initial data. This generalization derives from regularizing effects of Eikonal Hamiltonian (see \cite{b.g,b.q,b.r}). But our motivation is more in the case of equations (\ref{eq.i1})-(\ref{eq.i2}) where mutations are described by an integral operator. Then we can treat cases where the mutation rate $b(x,I)$ really depends on $x$, which was not available until now. The difficulty here is that Lipschitz bounds on the initial data are not propagated on $u_\epsilon$ and may blow up in finite time (see \cite{b.s,b.cc,b.hh} for regularity results for integral Hamiltonian). However, we achieve to control the Lipschitz norm by $-u_\epsilon$, that goes to infinity as $|x|$ goes to $+\infty$.\\

We do not discuss the uniqueness for equations (\ref{as.2}) and (\ref{as.i1}) in this paper. The latter is studied, for some particular cases, in \cite{b.a, b.h}.\\

A related, but different, situation arises in reaction-diffusion equations as in combustion (see \cite{b.m,b.k,b.ee,b.ff,b.t,b.gg}). A typical example is the Fisher-KPP equation, where the solution is a progressive front. The dynamics of the front is described by a level set of a solution of a Hamilton-Jacobi equation.\\

The paper is organized as follows. In section \ref{preresult} we state some existence results and bounds on $n_\epsilon$ and $I_\epsilon$. In section \ref{reg} we prove some regularity results for $u_\epsilon$ corresponding to equations~(\ref{eq.1})-(\ref{eq.2}). We show that $u_\epsilon$ are locally uniformly bounded and continuous. In section \ref{regint} we prove some analogous regularity results for $u_\epsilon$ corresponding to equations (\ref{eq.i1})-(\ref{eq.i2}). Finally, in section \ref{asymp} we describe the asymptotic behavior of $u_\epsilon$ and deduce the constrained Hamilton-Jacobi equation~(\ref{as.2})-(\ref{as.1.1}).

%---------------------------------------
\section{Preliminary results}
\label{preresult}

%---------------------------------------------

We recall the following existence results for $n_\epsilon$ and a priori bounds for $I_\epsilon$ (see also \cite{b.a, b.f}).

\begin{theorem}
\label{th.ex}
With the assumptions (\ref{eq.3})-(\ref{eq.5b}), and $I_m-C\epsilon^2 \leq I_\epsilon(0)\leq I_M+C\epsilon^2$, there is a unique solution $n_\epsilon \in \rC(\bR^+;L^1(\bR^d))$ to equations (\ref{eq.1})-(\ref{eq.2}) and it satisfies
\begin{equation}
\label{eq.7}
I_m'=I_m-C\epsilon^2 \leq I_\epsilon(t)\leq I_M+C\epsilon^2=I_M',
\end{equation}
where $C$ is a constant. This solution, $n_\epsilon(t,x)$, is nonnegative for all $t\geq 0$.
\end{theorem}

We recall a proof of this theorem in Appendix \ref{ap.1}. We have an analogous result for equations (\ref{eq.i1})-(\ref{eq.i2}):

\begin{theorem}
\label{th.exi}
With the assumptions (\ref{eq.5b}), (\ref{eq.5i})-(\ref{eq.6ij}), and $I_m\leq I_\epsilon(0)\leq I_M$, there is a unique solution $n_\epsilon \in \rC(\bR^+;L^1\cap L^\infty(\bR^d))$ to equations (\ref{eq.i1})-(\ref{eq.i2}) and it satisfies
\begin{equation}
\label{eq.7i}
I_m \leq I_\epsilon(t)\leq I_M.
\end{equation}
This solution, $n_\epsilon(t,x)$, is nonnegative for all $t\geq 0$.
\end{theorem}

This theorem can be proved with similar arguments as Theorem \ref{th.ex}. A uniform BV bound on $I_\epsilon(t)$ for equations (\ref{eq.1})-(\ref{eq.2}) is also proved in \cite{b.a}:

\begin{theorem}
\label{th.bv}
With the assumptions (\ref{eq.3})-(\ref{eq.6b}), we have additionally to the uniform bounds (\ref{eq.7}), the locally uniform BV and sub-Lipschitz bounds
\begin{equation}
\label{bv.1}
\frac{d}{dt}I_\epsilon(t)\geq -\epsilon \,C+e^{\frac{-Lt}{\epsilon}}\int \psi(x)n_\epsilon^0(x)\frac{R(x,I_\epsilon^0)}{\epsilon}dx,
\end{equation}
\begin{equation}
\label{bv.1b}
\frac{d}{dt}\varrho_\epsilon(t)\geq -Ct+\int (1+ \psi(x))n_\epsilon^0(x)\frac{R(x,I_\epsilon^0)}{\epsilon}dx,
\end{equation}
where $C$ and $L$ are positive constants and $\varrho_\epsilon(t)=\int_{\bR^d}n_\epsilon(t,x)dx$. Consequently, after extraction
of a subsequence, $I_\epsilon(t)$ converges a.e. to a function $I(t)$, as $\epsilon$ goes to $0$. The limit $I(t)$ is nondecreasing as soon as there exists a constant $C$ independent of $\epsilon$ such that

\begin{equation}
\label{bv.2}
\nonumber
\int \psi(x)n_\epsilon^0(x)\frac{R(x,I_\epsilon^0)}{\epsilon}\geq -Ce^{\frac{o(1)}{\epsilon}}.
\end{equation}

\end{theorem}

We also have a local BV bound on $I_\epsilon(t)$ for equations (\ref{eq.i1})-(\ref{eq.i2}):

\begin{theorem}
\label{th.bvi}
With the assumptions (\ref{eq.5b})-(\ref{eq.6ij}), we have additionally to the uniform bounds (\ref{eq.7i}), the locally uniform BV bound
\begin{equation}
\label{sL}
\frac{d}{dt}I_\epsilon(t)\geq -C'+e^{\frac{-L't}{\epsilon}}\int n_\epsilon^0(x)\frac{R(x,I_\epsilon^0)+b(x,I_\epsilon^0)}{\epsilon}dx,
\end{equation}
\begin{equation}
\label{sL2}
\int_0^T|\frac{d}{dt}I_\epsilon(t)|dt\leq 2C'T+C'',
\end{equation}
where $C'$, $C''$ and $L'$ are positive constants. Consequently, after extraction
of a subsequence, $I_\epsilon(t)$ converges a.e. to a function $I(t)$, as $\epsilon$ goes to $0$.

\end{theorem}

This theorem is proved in Appendix \ref{ap.2}.

%------------------------------------------------------
\section{Regularity results for equations (\ref{eq.1})-(\ref{eq.2})}
\label{reg}
%------------------------------------------------

In this section we study the regularity properties of $u_\epsilon=\epsilon \ln n_\epsilon$, where $n_\epsilon$ is the unique solution of equations (\ref{eq.1})-(\ref{eq.2}). We have

$$\partial_t n_\epsilon=\frac{1}{\epsilon}\partial_t u_\epsilon\, e^{\frac{u_\epsilon}{\epsilon}},\;\;\nabla n_\epsilon=\frac{1}{\epsilon}\nabla u_\epsilon\, e^{\frac{u_\epsilon}{\epsilon}},\;\;\triangle n_\epsilon=\big(\frac{1}{\epsilon}\triangle u_\epsilon+\frac{1}{\epsilon^2}|\nabla u_\epsilon|^2\big) e^{\frac{u_\epsilon}{\epsilon}}.$$
Consequently $u_\epsilon$ is a smooth solution to the following equation

\begin{equation}
\label{eq.H1}
\begin{cases}
\partial_t u_\epsilon-\epsilon\triangle u_\epsilon=|\nabla u_\epsilon|^2+R(x,I_\epsilon(t)),\;\;\;\;\; x\in \bR,\,t\geq0,\\
u_\epsilon(t=0)=\epsilon \ln n_\epsilon^0.
\end{cases}
\end{equation}

We have the following regularity results for $u_\epsilon$.

\begin{theorem}
\label{th1}
Assume (\ref{eq.3})-(\ref{as:n0-below}) and let $T>0$ be given. Set $D=B+(A^2+K_2)T$. Then we have $u_\epsilon\leq D^2$. For all $t_0>0$, $v_\epsilon=\sqrt{2D^2-u_\epsilon}$ are locally uniformly bounded and Lipschitz in $[t_0,T]\times \bR^d$,
\begin{equation}
\label{1}
|\nabla v_\epsilon| \leq C(T) (1+\frac{1}{\sqrt{t_0}}),
\end{equation}
where $C(T)$ is a constant depending on $T$, $K_1$, $K_2$, $A$ and $B$. Moreover, if we assume that $(u_\epsilon^0)_\epsilon:=\epsilon\ln(n_\epsilon^0)$ is a sequence of uniformly continuous functions, then $u_\epsilon$ are locally uniformly bounded and continuous in $[0,\infty[\times \bR^d$.
\end{theorem}

We prove Theorem \ref{th1} in several steps. We first prove an upper bound, then a regularizing effect in $x$, then local $L^\infty$ bounds, and finally a regularizing effect in $t$.

\subsection{An upper bound for $u_\epsilon$}
\label{sec.bo}

From assumption (\ref{eq.6b}) we have $u_\epsilon^0(x)\leq-A|x|+B$. We claim that, with $C= A^2+K_2$,

\begin{equation}
\label{eq.H3}
u_\epsilon(t,x)\leq -A|x|+B+Ct, \;\;\forall t\geq 0.
\end{equation}

Define $\phi(t,x)=-A|x|+B+Ct$. We have
$$\partial_t \phi-\epsilon\triangle \phi-|\nabla \phi|^2-R(x,I_\epsilon(t))\geq C+\epsilon\frac{A(d-1)}{|x|}-A^2-K_2\geq 0.$$

Here $K_2$ is an upper bound for $R(x,I)$ according to (\ref{eq.5b}). We have also $\phi(0,x)=-A|x|+B\geq u_\epsilon^0(x)$. So $\phi_\epsilon$ is a super-solution to (\ref{eq.H1}) and (\ref{eq.H3}) is proved.

\subsection{Regularizing effect in space}

\label{regx1}
Let $u=f(v)$, where $f$ is chosen later. We have
$$\partial_t u=f'(v)\partial_t v,\;\;
\partial_x u=f'(v)\partial_x v,\;\;\triangle u=f'(v)\triangle v+f''(v)|\nabla v|^2.$$

So equation (\ref{eq.H1}) becomes

\begin{eqnarray}
\label{eq.re2}
\partial_t v-\epsilon \triangle v-\left[\epsilon \frac{f''(v)}{f'(v)}+f'(v)\right] |\nabla v|^2=\frac{R(x,I)}{f'(v)}.
\end{eqnarray}

Define $p=\nabla v$. By differentiating (\ref{eq.re2}) we have

\begin{eqnarray}
\nonumber
\label{eq.re3}
\partial_t p_i-\epsilon \triangle p_i-2\left[\epsilon \frac{f''(v)}{f'(v)}+f'(v)\right] \nabla v\cdot \nabla p_i - \left[\epsilon \frac{f'''(v)}{f'(v)}-\epsilon \frac{f''(v)^2}{f'(v)^2}+f''(v)\right] |\nabla v|^2 p_i\\
\nonumber =-\frac{f''(v)}{f'(v)^2}R(x,I)p_i+\frac{1}{f'(v)}\frac{\partial R}{\partial x_i}.
\end{eqnarray}

We multiply the equation by $p_i$ and sum over $i$:

\begin{eqnarray}
\label{eq.re4}
\nonumber
\partial_t \frac{|p|^2}{2}-\epsilon \sum(\triangle p_i)p_i-2\left[\epsilon \frac{f''(v)}{f'(v)}+f'(v)\right] \nabla v\cdot \nabla \frac{|p|^2}{2} - \left[\epsilon \frac{f'''(v)}{f'(v)}-\epsilon \frac{f''(v)^2}{f'(v)^2}+f''(v)\right] |p|^4 \\
\nonumber =-\frac{f''(v)}{f'(v)^2}R(x,I)|p|^2+\frac{1}{f'(v)}\nabla_x R\cdot p.
\end{eqnarray}

First, we compute $\sum_i(\triangle p_i)p_i$.

\begin{alignat}{3}
\label{eq.re5}
\nonumber
\sum_i(\triangle p_i)p_i&=\sum_i\triangle \frac{p_i^2}{2}-\sum|\nabla p_i|^2\\
\nonumber
&=\triangle\frac{|p|^2}{2}-\sum|\nabla p_i|^2\\
\nonumber
&=|p|\triangle |p|+|\nabla|p||^2-\sum_i|\nabla p_i|^2.
\end{alignat}

We also have

\begin{alignat}{2}
\nonumber
|\nabla|p||^2=\sum_i \frac{|p\cdot \partial_{x_i}p|^2}{|p|^2}
\leq \sum_{i}|\partial_{x_i} p|^2=\sum_{i,j}|\partial_{x_i} p_j|^2=\sum_j|\nabla p_j|^2.
\end{alignat}

It follows that

\begin{eqnarray}
\nonumber
\label{eq.re6}
\sum_i(\triangle p_i)p_i\leq |p|\triangle|p|.
\end{eqnarray}

We deduce

\begin{eqnarray}
\label{eq.re7}
\partial_t |p|-\epsilon \triangle |p|-2\left[\epsilon \frac{f''(v)}{f'(v)}+f'(v)\right] p\cdot \nabla |p| - \left[\epsilon \frac{f'''(v)}{f'(v)}-\epsilon \frac{f''(v)^2}{f'(v)^2}+f''(v)\right] |p|^3 \\
\nonumber \leq-\frac{f''(v)}{f'(v)^2}R(x,I)|p|+\frac{1}{f'(v)}\nabla_x R\cdot \frac{p}{|p|}.
\end{eqnarray}

From (\ref{eq.H3}) we know that, for $0\leq t \leq T$, $u_\epsilon\leq D(T)^2$, where $D(T)=\sqrt{B+CT}$. Then we define $f(v)=-v^2+2D^2$, for $v$ positive, and thus

\begin{eqnarray}
\nonumber
\label{eq.re8}
D(T)\leq v,
\end{eqnarray}
\begin{eqnarray}
\label{eq.re9}
\nonumber
f'(v)=-2v, \;\;\;\;\text{and}\;\;\;\; |\frac{1}{f'(v)}|=\frac{1}{2v}\leq\frac{1}{2D},
\end{eqnarray}
\begin{eqnarray}
\label{eq.re10}
\nonumber
f''(v)=-2, \;\;\;\;\text{and}\;\;\;\;|\frac{f''(v)}{f'(v)^2}|=\frac{1}{2v^2}\leq\frac{1}{2D^2},
\end{eqnarray}
\begin{eqnarray}
\label{eq.re11}
f'''(v)=0, \;\;\;\;\nonumber-\left[\epsilon \frac{f'''(v)}{f'(v)}-\epsilon \frac{f''(v)^2}{f'(v)^2}+f''(v)\right]=2+\epsilon\frac{1}{v^2}>2.
\end{eqnarray}

From (\ref{eq.re7}), Theorem \ref{th.ex}, assumption (\ref{eq.5b}) and these calculations we deduce

\begin{eqnarray}
\nonumber
\label{eq.re13}
\frac{\partial |p|}{\partial t}-\epsilon \triangle |p|-2\left[\epsilon \frac{f''(v)}{f'(v)}+f'(v)\right] p\cdot \nabla |p| +2|p|^3 -\frac{K_2}{2D^2}|p|-\frac{K_2}{2D}\leq 0.
\end{eqnarray}

Thus for $\theta(T)$ large enough we can write

\begin{eqnarray}
\label{eq.re14}
\frac{\partial |p|}{\partial t}-\epsilon \triangle |p|-2\left[\epsilon \frac{f''(v)}{f'(v)}+f'(v)\right] p\cdot \nabla |p| +2(|p|-\theta)^3\leq 0.
\end{eqnarray}

Define the function

\begin{eqnarray}
\nonumber
\label{eq.re15}
y(t,x)=y(t)=\frac{1}{2\sqrt{ t}}+\theta.
\end{eqnarray}

Since $y$ is a solution to (\ref{eq.re14}), and $y(0)=\infty$ and $|p|$ being a sub-solution we have

\begin{eqnarray}
\nonumber
\label{eq.re16}
|p|(t,x) \leq y(t,x)=\frac{1}{2\sqrt{ t}}+\theta.
\end{eqnarray}

Thus for $v_\epsilon=\sqrt{2D^2-u_\epsilon}$, we have

\begin{equation}
\label{eq.re16b}
|\nabla v_\epsilon|(t,x) \leq \frac{1}{2\sqrt{t}}+\theta(T),\;\;\;\;0<t\leq T.
\end{equation}
See Appendix \ref{ap4} for more details on the comparison principle used above.

\subsection{Regularity in space of $u_\epsilon$ near $t=0$}

\label{regx2}
Assume that $u_\epsilon^0$ are uniformly continuous. We show that $u_\epsilon$ are uniformly continuous in space on $[0,T]\times \bR^d$.\\

For $\delta >0$ we prove that for $h$ small $|u_\epsilon(t,x+h)-u_\epsilon(t,x)|<\delta$. To do so define $w_\epsilon(t,x)=u_\epsilon(t,x+h)-u_\epsilon(t,x)$. Since $u_\epsilon^0$ are uniformly continuous, for $h$ small enough $|w_\epsilon(0,x)|<\frac{\delta}{2}$. Besides $w_\epsilon$ satisfies the following equation:

$$\partial_t w_\epsilon(t,x)-\epsilon \triangle w_\epsilon(t,x)-(\nabla u_\epsilon(t,x+h)+\nabla u_\epsilon(t,x))\cdot \nabla w_\epsilon(t,x)=R(x+h,I_\epsilon(t))-R(x,I_\epsilon(t)).$$

From Theorem \ref{th.ex} and using assumption (\ref{eq.5b}) we have

$$\partial_t w_\epsilon(t,x)-\epsilon \triangle w_\epsilon(t,x)-(\nabla u_\epsilon(t,x+h)+\nabla u_\epsilon(t,x))\cdot \nabla w_\epsilon(t,x)\leq K_2|h|.$$

Therefore by the maximum principle we arrive at

\begin{equation}
\nonumber
\label{eq.con}
\underset{\bR^d}{\max} |w_\epsilon(t,x)|<\underset{\bR^d}{\max}|w_\epsilon(0,x)|+K_2|h|t.
\end{equation}

So for $h$ small enough $|u_\epsilon(t,x+h)-u_\epsilon(t,x)|<\delta$ on $[0,T]\times \bR^d$.

\subsection{Local bounds for $u_\epsilon$}

\label{bound}
We show that $u_\epsilon$ are bounded on compact subsets of $ ]0,\infty[\times\bR^d$. We already know from section \ref{sec.bo} that $u_\epsilon$ is locally bounded from above. We show that it is also bounded from below on $\cC= [t_0,T]\times{\mathrm{B}}(0,R)$, for all $R>0$ and $0<t_0<T$.\\

From section \ref{sec.bo} we have $u_\epsilon(t,x)\leq -A|x|+B+CT$. So for $R$ large enough there exists $\epsilon_0$ such that for $\epsilon<\epsilon_0$
$$\int_{|x|>R}e^{\frac{u_\epsilon}{\epsilon}}dx<\int_{|x|>R}e^{\frac{-A|x|+B+CT}{\epsilon}}dx< \frac{I_m'}{2\psi_M}.$$

We have also from (\ref{eq.7}) that

$$\int_{\bR^d} e^{\frac{u_\epsilon}{\epsilon}}dx>\frac{I_m'}{\psi_M}.$$

We deduce that for $R$ large enough and for all $0<\epsilon<\epsilon_0$
$$\int_{|x|<R} e^{\frac{u_\epsilon}{\epsilon}}dx>\frac{I_m'}{2\psi_M}.$$

Therefore there exists $\epsilon_1>0$ such that, for all $\epsilon<\epsilon_1$
$$\exists x_0\in\bR^d;\;\;\;|x_0|<R, \; u_\epsilon(t,x_0)>-1, \;\;\text{thus}\;\;\;v_\epsilon(t,x_0)<\sqrt{2D^2+1}.$$

From Section \ref{regx1} we know that $v_\epsilon$ are locally uniformly Lipschitz
$$|v_\epsilon(t,x+h)-v_\epsilon(t,x)|<\big(C(T)+\frac{1}{2\sqrt{t_0}}\big)|h|,$$

Thus for all $(t,x)\in \cC$ and $\epsilon<\epsilon_1$

$$v_\epsilon(t,x)< E(t_0,T,R):= \sqrt{2D^2(T)+1}+2\big(C(T)+\frac{1}{2\sqrt{t_0}}\big)R.$$

It follows that

$$u_\epsilon(t,x)> 2D^2(T)-E^2(t_0,T,R).$$

We conclude that $u_\epsilon$ are uniformly bounded from below on $\cC$.\\

If we assume additionally that $u_\epsilon^0$ are uniformly continuous, with similar arguments we can show that $u_\epsilon$ are bounded on compact subsets of $[0,\infty[ \times\bR^d$. To prove the latter we use uniform continuity of $u_\epsilon$ instead of the Lipschitz bounds of $v_\epsilon$.

\subsection{Regularizing effect in time}
\label{sec.time}

From the above uniform bounds and continuity results we can also deduce uniform continuity in time i.e. for all $\eta>0$, there exists $\theta>0$ such that for all $(t,s,x)\in [0,T]\times [0,T] \times \mathrm{B}(0,\frac{R}{2})$, such that $0<t-s<\theta$, and for all $\epsilon<\epsilon_0$ we have:

$$|u_\epsilon(t,x)-u_\epsilon(s,x)|\leq 2\eta.$$

We prove this with the same method as that of Lemma $9.1$ in \cite{b.c} (see also \cite{b.kk} for another proof of this claim). We prove that for any $\eta>0$, we can find positive constants $A$, $B$ large enough such that, for any $x\in \mathrm{B}(0,\frac{R}{2})$, $s\in [0,T]$ and for every $\epsilon<\epsilon_0$,
\begin{equation}
\label{eq.regt}
u_\epsilon(t,y)-u_\epsilon(s,x)\leq \eta+A|x-y|^2+B(t-s),\;\;\;\text{for every}\; (t,y)\in [s,T]\times \mathrm{B}(0,{R}),
\end{equation}
and
\begin{equation}
\label{eq.regt2}
u_\epsilon(t,y)-u_\epsilon(s,x)\geq -\eta-A|x-y|^2-B(t-s),\;\;\;\text{for every}\; (t,y)\in [s,T]\times \mathrm{B}(0,{R}).
\end{equation}

We prove inequality (\ref{eq.regt}), the proof of (\ref{eq.regt2}) is analogous. We fix $(s,x)$ in $[0,T[\times\mathrm{B}(0,\frac{R}{2})$. Define 
$$\xi(t,y)=u_\epsilon(s,x)+\eta+A|y-x|^2+B(t-s),\;\;\;\;\; (t,y)\in [s,T[\times\mathrm{B}(0,R),$$
where $A$ and $B$ are constants to be determined. We prove that, for $A$ and $B$ large enough, $\xi$ is a super-solution to (\ref{eq.H1}) on $[s,T]\times \mathrm{B}(0,{R})$ and $\xi(t,y)>u_\epsilon(t,y)$ for $(t,y)\in \{s\}\times \mathrm{B}(0,{R})\cup [s,T]\times \partial\mathrm{B}(0,{R})$.\\

According to section \ref{bound}, $u_\epsilon$ are locally uniformly bounded, so we can take $A$ a constant such that for all $\epsilon<\epsilon_0$,

\begin{equation}
\nonumber
%\label{eq.regt3}
A\geq \frac{8\parallel u_\epsilon\parallel_{L^\infty([0,T]\times \mathrm{B}(0,{R}))}}{R^2}.
\end{equation}

With this choice, $\xi(t,y)>u_\epsilon(t,y)$ on $ [0,T]\times\partial \mathrm{B}(0,{R})$, for all $\eta$, $B$ and $x\in \mathrm{B}(0,\frac{R}{2})$. Next we prove that, for $A$ large enough, $\xi(s,y)>u_\epsilon(s,y)$ for all $y\in \mathrm{B}(0,{R})$. We argue by contradiction. Assume that there exists $\eta>0 $ such that for all constants $A$ there exists $y_{A,\epsilon} \in \mathrm{B}(0,{R})$ such that

\begin{equation}
\label{eq.regt4}
u_\epsilon(s,y_{A, \epsilon})-u_\epsilon(s,x)>\eta+A|y_{A, \epsilon}-x|^2.
\end{equation}

It follows that
\begin{equation}
\label{eq.regt5}
\nonumber
|y_{A, \epsilon}-x|\leq \sqrt{\frac{2M}{A}},
\end{equation}
where $M$ is a uniform upper bound for $\parallel u_\epsilon\parallel_{L^\infty([0,T]\times \mathrm{B}(0,{R}))}$. Now let $A \rightarrow \infty$. Then for all~$\epsilon$, $|y_{A, \epsilon}-x|\rightarrow 0$. According to Section \ref{regx2}, $u_\epsilon$ are uniformly continuous on space. Thus there exists $h>0$ such that if $|y_{A, \epsilon}-x|\leq h$ then $|u_\epsilon(s,y_{A, \epsilon})-u_\epsilon(s,x)|<\frac{\eta}{2}$, for all $\epsilon$. This is in contradiction with (\ref{eq.regt4}). Therefore $\xi(s,y)>u_\epsilon(s,y)$ for all $y\in \mathrm{B}(0,{R})$. Finally, noting that $R$ is bounded we deduce that for $B$ large enough, $\xi$ is a super-solution to (\ref{eq.H1}) in $[s,T]\times \mathrm{B}(0,{R})$. Since $u_\epsilon$ is a solution of (\ref{eq.H1}) we have

\begin{equation}
\nonumber
\label{eq.regt6}
u_\epsilon(t,y)\leq \xi(t,y)=u_\epsilon(s,x)+\eta+A|y-x|^2+B(t-s)\;\;\;\; \text{for all}\;\; (t,y)\in [s,T]\times \mathrm{B}(0,R).
\end{equation}

Thus (\ref{eq.regt}) is satisfied for $t\geq s$. We can prove (\ref{eq.regt2}) for $t\geq s$ analogously. Then we put $x=y$ and we conclude taking $\theta<\frac{\eta}{B}$.

%-------------------------------------
\section{Regularity results for equations (\ref{eq.i1})-(\ref{eq.i2})}
\label{regint}

%-----------------------------------

In this section we study the regularity properties of $u_\epsilon=\epsilon \ln n_\epsilon$, where $n_\epsilon$ is the unique solution of equations (\ref{eq.i1})-(\ref{eq.i2}) as given in Theorem \ref{th.exi}. From equation (\ref{eq.i1}) we deduce that $u_\epsilon$ is a solution to the following equation

\begin{equation}
\label{eq.iH1}
\begin{cases}
\partial_t u_\epsilon=R(x,I_\epsilon(t))+\int K(z) b(x+\epsilon z,I_\epsilon)e^{\frac{u_\epsilon(t,x+\epsilon z)-u_\epsilon(t,x)}{\epsilon}}dz,\;\;\;\;\; x\in \bR,\,t\geq0,\\
u_\epsilon(t=0)=\epsilon \ln n_\epsilon^0.
\end{cases}
\end{equation}

We have the following regularity results for $u_\epsilon$.
\begin{theorem}
\label{thi1}
 Let $n_\epsilon$ be the solution of (\ref{eq.i1})-(\ref{eq.i2}) with $n_\epsilon^0=e^{\frac{u\epsilon^0}{\epsilon}}$, and $u_\epsilon=\epsilon \ln(n_\epsilon)$. With the assumptions (\ref{eq.5b})-(\ref{eq.6ij}), and if we assume that $(u_\epsilon^0)_\epsilon$ is a sequence of uniformly bounded functions in $W^{1,\infty}$, then $u_\epsilon$ are locally uniformly bounded and Lipschitz in $[0,\infty[\times \bR^d$.

\end{theorem}

As in section \ref{reg} we prove Theorem \ref{thi1} in several steps. We first prove an upper and a lower bound on $u_\epsilon$, then local Lipschitz bounds in space and finally a regularity result in time.

\subsection{Upper and lower bounds on $u_\epsilon$}
\label{sec.ibo}

From assumption (\ref{eq.6b}) we have $u_\epsilon^0(x)\leq-A|x|+B$. As in section \ref{sec.bo} we claim that

\begin{equation}
\label{eq.iH3}
u_\epsilon(t,x)\leq -A|x|+B+Ct, \;\;\forall t\geq 0.
\end{equation}

Define $v(t,x)=-A|x|+B+Ct$, where $C= b_M\int K(z)e^{A|z|}dz+K_2$. Using (\ref{eq.5b}) and (\ref{eq.6i}) we have

$$\partial_t v-R(x,I_\epsilon(t))-\int K(z) b(x+\epsilon z,I_\epsilon)e^{\frac{v(t,x+\epsilon z)-v(t,x)}{\epsilon}}dz\geq C-K_2-b_M\int K(z)e^{A|z|}dz\geq 0.$$

We also have $v(0,x)=-A|x|+B\geq u_\epsilon^0(x)$. So $v$ is a supersolution to (\ref{eq.iH1}). Since (\ref{eq.i1}) verifies the comparison property, equation (\ref{eq.iH1}) verifies also the comparison property, i.e. if $v$ and $u$ are respectively super and subsolutions of (\ref{eq.iH1}) then $u\leq v$. Thus (\ref{eq.iH3}) is proved.\\

To prove a lower bound on $u_\epsilon $ we assume that $u_\epsilon^0$ are locally uniformly bounded. Then from equation (\ref{eq.iH1}) and assumption (\ref{eq.5b}) we deduce
$$\partial_t u_\epsilon(t,x)\geq -K_2, $$
and thus 
\begin{eqnarray}
\nonumber u_\epsilon(t,x)\geq-\|u_\epsilon^0\|_{L^\infty({\mathrm{B}}(0,R))}-K_2t, \;\;\;\; \forall x\in {\mathrm{B}}(0,R).
\end{eqnarray}

Moreover, $|\nabla u_\epsilon^0|$ being bounded, we can give a lower bound in $\bR^d$

\begin{eqnarray}
\label{lowerb}
u_\epsilon(t,x)\geq \underset{\epsilon}\inf \,u_\epsilon^0(0)-\|\nabla u_\epsilon^0\|_{L^\infty}|x|-K_2t, \;\;\;\; \forall x\in \bR^d.
\end{eqnarray}

\subsection{Lipschitz bounds}
\label{lipx}

Here we assume that $u_\epsilon$ is differentiable in $x$ (See \cite{b.hh}). See also Appendix \ref{ap3} for a proof without any regularity assumptions on $u_\epsilon$.\\

Let $p_\epsilon=\nabla u_\epsilon \cdot \chi$, where $\chi$ is a fixed unit vector. By differentiating (\ref{eq.iH1}) with respect to $\chi$ we obtain

\begin{align}
\nonumber
\partial_t p_\epsilon(t,x)&=\nabla R(x,I_\epsilon(t))\cdot \chi+\int K(z) \nabla b(x+\epsilon z,I_\epsilon)\cdot \chi \,e^{\frac{u_\epsilon(t,x+\epsilon z)-u_\epsilon(t,x)}{\epsilon}}dz\\
\nonumber
&+\int K(z)b(x+\epsilon z,I_\epsilon)\frac{p_\epsilon(t,x+\epsilon z)-p_\epsilon(t,x)}{\epsilon}e^{\frac{u_\epsilon(t,x+\epsilon z)-u_\epsilon(t,x)}{\epsilon}}dz.
\end{align}

Thus, using assumptions (\ref{eq.5b}) and (\ref{eq.6i}), we have
\begin{align}
\label{bornp}
\partial_t p_\epsilon(t,x)&\leq K_2+L_1\int K(z) b(x+\epsilon z,I_\epsilon)e^{\frac{u_\epsilon(t,x+\epsilon z)-u_\epsilon(t,x)}{\epsilon}}dz\\
\nonumber
&+\int K(z)b(x+\epsilon z,I_\epsilon)\frac{p_\epsilon(t,x+\epsilon z)-p_\epsilon(t,x)}{\epsilon}e^{\frac{u_\epsilon(t,x+\epsilon z)-u_\epsilon(t,x)}{\epsilon}}dz.
\end{align}

Define $w_\epsilon(t,x)=p_\epsilon(t,x)+L_1u_\epsilon(t,x)$ and $\Delta_\epsilon(t,x,z)=\frac{u_\epsilon(t,x+\epsilon z)-u_\epsilon(t,x)}{\epsilon}$. From (\ref{bornp}) and (\ref{eq.iH1}) we deduce 
\begin{align}
\nonumber
&\partial_t w_\epsilon- K_2(1+L_1)-\int K(z)b(x+\epsilon z,I_\epsilon)\frac{w_\epsilon(t,x+\epsilon z)-w_\epsilon(t,x)}{\epsilon}e^{\Delta_\epsilon(t,x,z)}dz\\
\nonumber
&\leq 2L_1\int K(z) b(x+\epsilon z,I_\epsilon)e^{\Delta_\epsilon(t,x,z)}dz\\
\nonumber
&-L_1\int K(z)b(x+\epsilon z,I_\epsilon)\Delta_\epsilon(t,x,z)e^{\Delta_\epsilon(t,x,z)}dz\\
\nonumber
&= L_1\int K(z) b(x+\epsilon z,I_\epsilon)e^{\Delta_\epsilon(t,x,z)}\big(2-\Delta_\epsilon(t,x,z)\big)dz\\
\nonumber 
&\leq L_1 b_M e,
\end{align}
noticing that $e$ is the maximum of the function $g(t)=e^t(2-t)$ in $\bR$. Therefore by the maximum principle, with $C_1=K_2(1+L_1)+ L_1 b_M e$, we have

$$w_\epsilon(t,x)\leq C_1t+\underset{\bR^d}{\max}\;w_\epsilon(0,x).$$

It follows that
\begin{align}
\label{local} p_\epsilon(t,x)&\leq C_1t+\parallel \nabla u_\epsilon^0\parallel_{L^\infty}+L_1(B+Ct)+L_1\big( \|\nabla u_\epsilon^0\|_{L^\infty}|x|+K_2t-u_\epsilon^0(x=0) \big)\\
\nonumber
&= C_2t+C_3|x|+C_4,
\end{align}
where $C_2$, $C_3$ and $C_4$ are constants. Since this bound is true for any $|\chi| = 1$, we obtain a local bound on $|\nabla u_\epsilon|$.

\subsection{Regularity in time}
In section \ref{lipx} we proved that $u_\epsilon$ is locally uniformly Lipschitz in space. From this we can deduce that $\partial_t u_\epsilon$ is also locally uniformly bounded.\\

Let $\mathcal{C}=[0,T]\times \mathrm{B}(x_0,R)$ and $S_1$ be a constant such that $\parallel u_\epsilon \parallel_{L^\infty(\mathcal{C})}<S_1$ for all $\epsilon>0$. Assume that $R'$ is a constant large enough such that we have $u_\epsilon(t,x)<-S_1$ in $[0,T]\times \bR^d\backslash\mathrm{B}(x_0,R')$. According to (\ref{eq.iH3}) there exists such constant $R'$. We choose a constant $S_2$ such that $\parallel \nabla u_\epsilon \parallel_{L^\infty([0,T]\times \mathrm{B}(x_0,R'))}<S_2$ for all $\epsilon>0$. We deduce

\begin{align}
\nonumber
|\partial_t u_\epsilon|&\leq |R(x,I_\epsilon(t))|+\int K(z) b(x+\epsilon z,I_\epsilon)e^{\frac{u_\epsilon(t,x+\epsilon z)-u_\epsilon(t,x)}{\epsilon}}\big( \mathds{1}_{|x+\epsilon z|<R'}+\mathds{1}_{|x+\epsilon z|\geq R'}\big)dz\\
\nonumber
&\leq K_2+b_M \int K(z)e^{S_2 |z|} \mathds{1}_{|x+\epsilon z|<R'}dz+b_M \int K(z)\mathds{1}_{|x+\epsilon z|\geq R'}dz\\
\nonumber
&\leq K_2+b_M\big(1+\int K(z)e^{S_2|z|}dz\big).
\end{align}

This completes the proof of Theorem \ref{thi1}.

%------------------------------------------------

\section{Asymptotic behavior of $u_\epsilon$}
\label{asymp}

%------------------------------------------------

Using the regularity results in sections \ref{reg} and \ref{regint}, we can now describe the asymptotic behavior of $u_\epsilon$ and prove Theorems \ref{th.3} and \ref{th.i2}. Here we prove Theorem \ref{th.3}. The proof of Theorem \ref{th.i2} is analogous, except the limit of the integral term in equation (\ref{as.i1}). The latter has been studied in \cite{b.e,b.s,b.h,b.l}.

\begin{proof}[Proof of theorem \ref{th.3}]

\textbf{step 1 (Limit)} According to section \ref{reg}, $u_\epsilon$ are locally uniformly bounded and continuous. So by Arzela-Ascoli Theorem after extraction of a subsequence, $u_\epsilon$ converges locally uniformly to a continuous function $u$.\\

\textbf{step 2 (Initial condition)} We proved that if $u_\epsilon^0$ are uniformly continuous then $u_\epsilon$ will be locally uniformly bounded and continuous in $[0,T]\times \bR^d$. Thus we can apply Arzela-Ascoli near $t=0$ as well. Therefore we have $u(0,x)=\underset{\epsilon\rightarrow0} {\lim}\; u_\epsilon(0,x)=u^0(x)$.\\

\textbf{step 3 \big($\underset{x\in \bR^d}\max \,u=0$\big)} Assume that for some $t,x$ we have $0<a\leq u(t,x)$. Since $u$ is continuous $u(t,y)\geq \frac{a}{2}$ on $\mathrm{B}(x,r)$, for some $r>0$. Thus we have $n_\epsilon(t,y)\rightarrow \infty$, while $\epsilon\rightarrow 0$. Therefore $I_\epsilon(t)\rightarrow \infty$ while $\epsilon\rightarrow 0$. This is a contradiction with (\ref{eq.7}).\\

To prove that $\underset{x\in \bR^d}{\max} \;u(t,x)=0$, it suffices to show that $\underset{\epsilon\rightarrow0} {\lim}\; n_\epsilon(t,x)\neq0$, for some $x\in \bR^d$. From (\ref{eq.H3}) we have
$$u_\epsilon(t,x)\leq -A|x|+B+Ct.$$

It follows that for $M$ large enough

\begin{eqnarray}
\label{aaa}\underset{\epsilon\rightarrow0} {\lim}\int_{|x|>M}n_\epsilon(t,x)dx\leq \underset{\epsilon\rightarrow0} {\lim}\int_{|x|>M}e^{\frac{-A|x|+B+Ct}{\epsilon}}=0.
\end{eqnarray}

From this and (\ref{eq.7}) we deduce $$\underset{\epsilon\rightarrow0} {\lim}\int_{|x|\leq M}n_\epsilon(t,x)dx\geq\frac{I'_m}{\psi_M}.$$

If $u(t,x)<0$ for all $|x|<M$ then $\underset{\epsilon\rightarrow0}{\lim}\;e^{\frac{u_\epsilon(t,x)}{\epsilon}}=0$ and thus $\underset{\epsilon\rightarrow0}{\lim}\int_{|x|\leq M}n_\epsilon(t,x)dx=0$. This is a contradiction with (\ref{aaa}).
It follows that $\underset{x\in \bR^d}{\max} \;u(t,x)=0,\;\;\;\forall t>0$.\\

\textbf{step 4 ($supp\; n(t,\cdot)\subset \{u(t,\cdot)=0\}$)} Assume that $u(t_0,x_0)=-a<0$. Since $u_\epsilon$ are uniformly continuous in a small neighborhood of $(t_0,x_0)$, $(t,x)\in [t_0-\delta,t_0+\delta]\times \mathrm{B}(x_0,\delta)$, we have $u_\epsilon (t,x)\leq-\frac{a}{2}<0$ for $\epsilon$ small. We deduce that $\int_{ [t_0-\delta,t_0+\delta]\times \mathrm{B}(x_0,\delta)}n \,dt dx=\int_{ [t_0-\delta,t_0+\delta]\times \mathrm{B}(x_0,\delta)}\underset{\epsilon \rightarrow 0 }{\lim } \;e^{\frac{u_\epsilon(t,x)}{\epsilon}}dt dx=0$. Therefore we have $supp\; n(t,\cdot)\subset \{u(t,\cdot)=0\}$ for almost every $t$. \\

\textbf{step 5 (Limit equation)} Finally we recall, following \cite{b.a}, how to pass to the limit in the equation. Since $u_\epsilon$ is a solution to (\ref{eq.H1}), it follows that $\phi_\epsilon(t,x)=u_\epsilon(t,x)-\int_0^t R(x,I_\epsilon(s))ds$ is a solution to the following equation
\begin{eqnarray}
\nonumber
\partial_t \phi_\epsilon(t,x)-\epsilon\triangle \phi_\epsilon(t,x)-|\nabla \phi_\epsilon(t,x)|^2-2\nabla \phi_\epsilon(t,x)\cdot \int_0^t \nabla R(x,I_\epsilon(s))ds\\
\nonumber
= \epsilon \int_0^t \triangle R(x,I_\epsilon(s))ds+|\int_0^t \nabla R(x,I_\epsilon(s))ds|^2.
\end{eqnarray}

Note that we have $I_\epsilon(s)\rightarrow I(s)$ for all $s\geq 0$ as $\epsilon$ goes to $0$, and on the other hand, the function $R(x,I)$ is smooth. It follows that we have the locally uniform limits
$$\underset{\epsilon\rightarrow0}{\lim}\int_0^t R(x,I_\epsilon(s))ds= \int_0^t R(x,I(s))ds,$$
$$\underset{\epsilon\rightarrow0}{\lim}\int_0^t \nabla R(x,I_\epsilon(s))ds=\int_0^t \nabla R(x,I(s))ds,$$
$$\underset{\epsilon\rightarrow0}{\lim}\int_0^t \triangle R(x,I_\epsilon(s))ds=\int_0^t \triangle R(x,I(s))ds,$$
for all $t\geq 0$. Moreover the functions $\int_0^t R(x,I(s))ds$, $\int_0^t \nabla R(x,I(s))ds$ and $\int_0^t \triangle R(x,I(s))ds$ are continuous. According to step 1, $u_\epsilon(t,x)$ converge locally uniformly to the continuous function $u(t,x)$ as $\epsilon$ vanishes. Therefore $\phi_\epsilon(t,x)$ converge locally uniformly to the continuous function $\phi(t,x)=u(t,x)-\int_0^t R(x,I(s))ds$ as $\epsilon$ vanishes. It follows that $\phi(t,x)$ is a viscosity solution to the equation

\begin{eqnarray}
\nonumber
\partial_t \phi(t,x)-|\nabla \phi(t,x)|^2-2\nabla \phi(t,x)\cdot \int_0^t \nabla R(x,I(s))ds\\
\nonumber
= |\int_0^t \nabla R(x,I)ds|^2.
\end{eqnarray}

In other words $u(t,x)$ is a viscosity solution to the following equation

\begin{eqnarray}
\nonumber
\partial_t u(t,x)=|\nabla u(t,x)|^2+R(x,I(t)).
\end{eqnarray}

\end{proof}

%-----------------------------------------
\appendix
%-----------------------------------------

%-----------------------------------------
\section{Proof of theorem \ref{th.ex}}
\label{ap.1}
%-----------------------------------------

\subsection{Existence}
Let $T>0$ be given and $\mathrm{A}$ be the following closed subset:
$$\mathrm{A}=\{u\in \rC\big([0,T],L^1(\bR^d)\big),\; u\geq 0,\; \parallel u(t,\cdot)\parallel_{L^1}\leq a\},$$
where $a= \left( \int n_\epsilon^0 dx\right) e^{\frac{K_2T}{\epsilon}}$. Let $\Phi$ be the following application:

\begin{equation}
\nonumber\label{eq.8}
\Phi:\mathrm{A}\rightarrow\mathrm{A}
\end{equation}
\begin{equation}
\nonumber u\mapsto v,
\end{equation}
where $v$ is the solution to the following equation
\begin{equation}
\label{eq.9}
\begin{cases}
\partial_t v-\epsilon \triangle v=\frac{v}{\epsilon}\bar{R}(x,I_u(t)),\;\;\;\;\; x\in \bR,\,t\geq0,\\
v(t=0)=n_\epsilon^0.
\end{cases}
\end{equation}
\begin{equation}
\label{eq.10}
I_u(t)=\int_{\bR^d}\psi(x)u(t,x)dx,
\end{equation}
and $\bar{R}$ is defined as below
\begin{equation}
\nonumber
\bar{R}(x,I)=
\begin{cases}
R(x,I)\;\;\;\; &\text{if} \;\;\; \frac{I_m}{2}<I<2I_M,\\
R(x,2I_M)\;\;\;\; &\text{if} \;\;\; 2I_M\leq I,\\
R(x,\frac{I_m}{2})\;\;\;\; &\text{if} \;\;\; I\leq \frac{I_m}{2}.
\end{cases}
\end{equation}

We prove that
\begin{enumerate}
\item \label{i.a}$\Phi$ defines a mapping of $A$ into itself,
\item \label{i.b}$\Phi$ is a contraction for $T$ small.
\end{enumerate}

With these properties, we can apply the Banach-Picard fixed point theorem and iterate the construction with $T$ fixed.\\

Assume that $u\in \mathrm{A}$. In order to prove (a) we show that $v$, the solution to (\ref{eq.9}), belongs to $\mathrm{A}$. By the maximum principle we know that $v\geq 0$. To prove the $L^1$ bound we integrate (\ref{eq.9})
$$\frac{d}{dt}\int v dx=\int\frac{v}{\epsilon}\bar{R}(x,I_u(t))dx\leq \frac{1}{\epsilon}\,\underset{x\in \bR^d}{\max}\, \bar{R}(x,I_u(t))\int v dx\leq \frac{K_2}{\epsilon}\int v dx,$$
and we conclude from the Gronwall Lemma that

$$\parallel v \parallel_{L^1}\leq \left( \int n_\epsilon^0 dx\right) e^{\frac{K_2T}{\epsilon}}=a.$$

Thus (a) is proved. It remains to prove (b). Let $u_1,\; u_2\in \mathrm{A}$, $v_1=\Phi(u_1)$ and $v_2=\Phi(u_2)$. We have

\begin{equation}
\nonumber
\partial_t(v_1-v_2)-\epsilon \triangle (v_1-v_2)=\frac{1}{\epsilon}\left[ (v_1-v_2)\bar{R}(x,I_{u_1})+v_2\left( \bar{R}(x,I_{u_1})-\bar{R}(x,I_{u_2})\right) \right].
\end{equation}

Noting that $\parallel v_2 \parallel_{L^1}\leq a$, and $| \bar{R}(x,I_{u_1}) -\bar{R}(x,I_{u_2})| \leq K_1|I_{u_1}-I_{u_2}|\leq K_1\psi_M \parallel u_1-u_2 \parallel_{L^1}$ we obtain

$$\frac{d}{dt}\parallel v_1-v_2\parallel_{L^1}\leq \frac{K_2}{\epsilon}\parallel v_1-v_2\parallel_{L^1}+\frac{a K_1\psi_M}{\epsilon}\parallel u_1-u_2\parallel_{L^1}.$$

Using $v_1(0,\cdot)=v_2(0,\cdot)$ we deduce

$$\parallel v_1-v_2\parallel_{L^\infty_t L^1_x}\leq \frac{aK_1\psi_M}{K_2}(e^{\frac{K_2T}{\epsilon}}-1)\parallel u_1-u_2\parallel_{L^\infty_t L^1_x}.$$

Thus, for $T$ small enough such that $e^{\frac{K_2T}{\epsilon}}(e^{\frac{K_2T}{\epsilon}}-1)<\frac{K_2}{2K_1\psi_M\int n_\epsilon^0}$, $\Phi$ is a contraction. Therefore $\Phi$ has a fixed point and there exists $n_\epsilon\in \mathrm{A}$ a solution to the following equation

\begin{equation}
\nonumber
\begin{cases}
\partial_t n_\epsilon-\epsilon \triangle n_\epsilon=\frac{n_\epsilon}{\epsilon}\bar{R}(x,I(t)),\;\;\;\;\; x\in \bR,\,0\leq t\leq T,\\
n_\epsilon(t=0)=n_\epsilon^0.
\end{cases}
\end{equation}

\begin{equation}
\nonumber
I(t)=\int_{\bR^d}\psi(x)n_\epsilon(t,x)dx,
\end{equation}

With the same arguments as \ref{ap1.2} we prove that $\frac{I_m}{2}<I(t)<2I_M$ and thus $n_\epsilon$ is a solution to equations (\ref{eq.1})-(\ref{eq.2}) for $t\in[0,T]$. We fix $T$ small enough such that $e^{\frac{K_2T}{\epsilon}}(e^{\frac{K_2T}{\epsilon}}-1)<\frac{K_2\psi_m}{4K_1\psi_M I_M}$. Then we can iterate in time and find a global solution to equations (\ref{eq.1})-(\ref{eq.2}).\\

\subsection{Uniform bounds on $I_\epsilon(t)$}
\label{ap1.2}

We have

\begin{align}
\nonumber
\frac{d I_\epsilon}{dt}=\frac{d}{dt} \int_{\bR^d}\psi(x)n_\epsilon(t,x)dx=
\epsilon\int_{\bR^d}\psi(x)\triangle n_\epsilon(t,x) dx+ \frac{1}{\epsilon}\int_{\bR^d}\psi(x)n_\epsilon(t,x) R(x,I_\epsilon(t))dx.
\end{align}

We define $\psi_L=\chi_L\cdot \psi \in \mathbf{W}_{2,c}^\infty(\bR^d)$, where $\chi_L$ is a smooth function with a compact support such that $\chi_L|_{\mathrm{B}(0,L)}\equiv 1$, $\chi_L|_{\mathbb{R}\backslash\mathrm{B}(0,2L)}\equiv 0$. Then by integration by parts we find

\begin{align}
\nonumber
%\label{A2}
\int_{\bR^d}\psi_L(x) \triangle n_\epsilon(t,x) dx=\int_{\bR^d}\triangle\psi_L(x) n_\epsilon(t,x) dx.
\end{align}

As $L\rightarrow \infty$, $\psi_L$ converges to $\psi$ in $W^{2,\infty}_{\text{loc}}(\bR^d)$. Therefore we obtain

\begin{align}
\nonumber%\label{A3}
\underset{L\rightarrow \infty}{\lim}\int_{\bR^d}\triangle\psi_L(x) n_\epsilon dx &=\int_{\bR^d}\triangle\psi(x) n_\epsilon dx,\\
\nonumber
\underset{L\rightarrow \infty}{\lim}\int_{\bR^d}\psi_L(x) \triangle n_\epsilon(t,x) dx &=\int_{\bR^d}\psi(x) \triangle n_\epsilon(t,x) dx.
\end{align}

From these calculations we conclude

\begin{align}
\nonumber
%\label{A4}
\frac{d I_\epsilon}{dt}=
\epsilon\int_{\bR^d}\triangle\psi(x) n_\epsilon(t,x) dx+ \frac{1}{\epsilon}\int_{\bR^d}\psi(x)n_\epsilon(t,x) R(x,I_\epsilon(t))dx.
\end{align}

It follows that
\begin{equation}
\nonumber%\label{A5}
-\epsilon \frac{C_1}{\psi_m} I_\epsilon+\frac{1}{\epsilon}I_\epsilon\min_{x\in \bR^d} {R(x,I_\epsilon)}\leq \frac{dI_\epsilon }{dt}\leq\epsilon \frac{C_1}{\psi_m} I_\epsilon+\frac{1}{\epsilon}I_\epsilon\max_{x\in \bR^d} {R(x,I_\epsilon)}.
\end{equation}

Let $C= \frac{ C_1 K_1}{\psi_m}$. As soon as $I_\epsilon$ overpasses $I_M+C\epsilon^2$, we have $R(x,I_\epsilon)<-\frac{C\epsilon^2}{K_1}=-\epsilon^2\frac{C_1}{\psi_m}$ and thus $\frac{dI_\epsilon }{dt}$ becomes negative. Similarly, as soon as $I_\epsilon$ becomes less than $I_m-C\epsilon^2$, $\frac{dI_\epsilon }{dt}$ becomes positive. Thus (\ref{eq.7}) is proved.

%------------------------------------

\section{A locally uniform BV bound on $I_\epsilon$ for equations (\ref{eq.i1})-(\ref{eq.i2})}
\label{ap.2}

%------------------------------------

In this appendix we prove Theorem \ref{th.bvi}. We first integrate (\ref{eq.i1}) over $\bR^d$ to obtain

\begin{equation}
\nonumber
\frac{d}{dt} I_\epsilon(t)=\frac{1}{\epsilon}\int n_\epsilon(t,x)\big(R\left(x,I_\epsilon(t)\right)+b\left( x,I_\epsilon(t)\right) \big)dx.
\end{equation}

Define $J_\epsilon(t)=\frac{d}{dt} I_\epsilon(t)$. We differentiate $J_\epsilon$ and we obtain

\begin{align}
\nonumber
\frac{d}{dt} J_\epsilon(t)&=\frac{1}{\epsilon}J_\epsilon(t)\int n_\epsilon(t,x) \frac{\partial (R+b)}{\partial I}(x,I_\epsilon(t)) dx \\
\nonumber
&+\frac{1}{\epsilon^2}\int \big(R(x,I_\epsilon)+b(x,I_\epsilon)\big)\big[n_\epsilon(t,x)R(x,I_\epsilon)+\int K_\epsilon(y-x)b(y,I_\epsilon)n_\epsilon(t,y)dy \big]dx.
\end{align}

We rewrite this equality in the following form

\begin{align}
\nonumber
\frac{d}{dt} J_\epsilon(t)&=\frac{1}{\epsilon}J_\epsilon(t)\int n_\epsilon(t,x) \frac{\partial (R+b)}{\partial I}\big(x,I_\epsilon(t)\big) dx +\frac{1}{\epsilon^2}\int n_\epsilon(t,x)\big(R\big(x,I_\epsilon(t)\big)+b\big(x,I_\epsilon(t)\big)\big)^2 dx\\
\nonumber
&+\frac{1}{\epsilon^2}\int \int K_\epsilon(y-x)\big(R\big(x,I_\epsilon(t)\big)-R\big(y,I_\epsilon(t)\big)\big)b\big(y,I_\epsilon(t)\big)n_\epsilon(t,y)dy dx\\
\nonumber
&+\frac{1}{\epsilon^2}\int \int K_\epsilon(y-x)\big(b\big(x,I_\epsilon(t)\big)-b\big(y,I_\epsilon(t)\big)\big)b\big(y,I_\epsilon(t)\big)n_\epsilon(t,y)dy dx.
\end{align}

It follows that
\begin{align}
\nonumber
\frac{d}{dt} J_\epsilon(t)&\geq \frac{1}{\epsilon}J_\epsilon(t)\int n_\epsilon(t,x) \frac{\partial (R+b)}{\partial I}\big(x,I_\epsilon(t)\big) dx +\frac{1}{\epsilon^2}\int n_\epsilon(t,x)\big(R\big(x,I_\epsilon(t)\big)+b\big(x,I_\epsilon(t)\big)\big)^2 dx\\
\nonumber
&-\frac{K_2+b_M \,L_1}{\epsilon}\int \int K(z) |z| b\big(x+\epsilon z,I_\epsilon(t)\big)n_\epsilon(t,x+\epsilon z)dz dx\\
% \nonumber
% &-\frac{b_M \,L_1}{\epsilon}\int \int K(z) |z| b(x+\epsilon z,I_\epsilon)n_\epsilon(t,x+\epsilon z)dz dx\\
\nonumber
&\geq \frac{1}{\epsilon}J_\epsilon(t)\int n_\epsilon(t,x) \frac{\partial (R+b)}{\partial I}\big(x,I_\epsilon(t)\big) dx +\frac{1}{\epsilon^2}\int n_\epsilon(t,x)\big(R\big(x,I_\epsilon(t)\big)+b\big(x,I_\epsilon(t)\big)\big)^2 dx-\frac{C_1}{\epsilon},
\end{align}
where $C_1$ is a positive constant. Consequently, using (\ref{eq.6ij}) we obtain

$$\frac{d}{dt}(J_\epsilon(t))_- \leq \frac{C_1}{\epsilon}-\frac{C_2}{\epsilon}(J_\epsilon(t))_-,$$

with $(J_\epsilon(t))_-=\max(0,-J_\epsilon(t))$. From this inequality we deduce

$$(J_\epsilon(t))_-\leq \frac{C_1}{C_2}+(J_\epsilon(0))_-e^{-\frac{C_2t}{\epsilon}}.$$

With similar arguments we obtain

$$(J_\epsilon(t))_+\geq -\frac{C'_1}{C'_2}+(J_\epsilon(0))_+e^{-\frac{C'_2t}{\epsilon}},$$

with $(J_\epsilon(t))_+=\max(0,J_\epsilon(t))$. Thus (\ref{sL}) is proved. Finally, we deduce the locally uniform BV bound (\ref{sL2})

\begin{align}
\nonumber
\int_0^T |\frac{d}{dt}I_\epsilon(t)| dt&=\int_0^T \frac{d}{dt}I_\epsilon(t) dt+2\int_0^T (\frac{d}{dt}I_\epsilon(t))_- dt\\
\nonumber
&\leq I_M-I_m+2C'T+O(1).
\end{align}

%--------------------------------------

\section{Complement to the proof of the regularizing effect   (\ref{1})}
\label{ap4}
%--------------------------------------

In this section, we provide some details for the comparison principle used in the proof of (\ref{1}). In Subsection \ref{regx1} we proved that $p=\nabla v$ satisfies the following (see the inequality (\ref{eq.re14}))
$$
\frac{\partial |p|}{\partial t}-\epsilon \triangle |p|-2\left[\frac{\epsilon}{v} -2v \right] p\cdot \nabla |p| +2(|p|-\theta)^3\leq 0.
$$
To apply the comparison principle we first claim the following lemma that we will prove at the end of this section.
\begin{lemma}
\label{lem:lb}
Assume  (\ref{eq.5b}) and (\ref{as:n0-below}). Then, there exist positive constants $A_1$, $B_1$ and $D_1$ such that, for all $t_1>0$ and $\e\leq 1$,
\begin{equation}
\label{low-b-u}
-\frac{A_1|x|^2+B_1+D_1t}{t_1} \leq u_\e(t,x),\qquad \text{for $(t,x)\in (t_1,+\infty)\times \bR^d$}.
\end{equation}
\end{lemma}
The above lemma implies that
$$
D(T)\leq v_\e \leq \sqrt{ 2D^2+\frac{1}{t_1}(B_1+D_1T+A_1|x|^2)},\qquad \text{for $(t,x)\in (t_1,+\infty)\times \bR^d$}.
$$
We deduce that, for some positive constants $A_2$ and $D_2(T)$,
\begin{equation}
\label{eq:p}
\frac{\partial |p|}{\partial t}-\epsilon \triangle |p|- \frac{1}{\sqrt{t_1}}\left[A_2|x|+D_2(T) \right] |p|\cdot \nabla |p| +2(|p|-\theta)^3\leq 0,\qquad \text{for $(t,x)\in (t_1,+\infty)\times \bR^d$}.\end{equation}
Define, for $(t,x)\in (t_1,T]\times B_R(0)$ and for $A_3$ to be chosen later,
 $$
z(t,x)=\f{1}{2\sqrt{t-t_1}}+\f{A_3R^2}{\sqrt{t_1}(R^2-|x|^2)}+\theta.
 $$
We prove that,
for $A_3=A_3(T)$ chosen large enough, $z$ is a strict supersolution of (\ref{eq:p}) in $(t_1,T]\times B_R(0)$.
 To this end, we compute
 $$
 \p_{t} z(t,x)=-\f{1}{4(t-t_1)\sqrt{t-t_1}},
 $$
 $$
 \nabla z(t,x)=\f{2A_3R^2 x}{\sqrt{t_1}(R^2-|x|^2)^2},
 $$
  $$
 \Delta z(t,x)=\f{2A_3R^2 }{\sqrt{t_1}(R^2-|x|^2)^2} +\f{8A_3R^2 |x|^2}{\sqrt{t_1}(R^2-|x|^2)^3}.
 $$
We then replace this in (\ref{eq:p}) to obtain
$$
\begin{array}{rl}
 & \frac{\partial z}{\partial t}-\e \Delta z-\frac{1}{\sqrt{t_1}}( A_2|x|+D_2(T)) |z \nabla z| +2(z-\theta)^3 =\\
  &-\f{1}{4(t-t_1)\sqrt{t-t_1}}-\e\big( \f{2A_3R^2}{\sqrt{t_1}(R^2-|x|^2)^2}+\f{8A_3R^2|x|^2}{\sqrt{t_1}(R^2-|x|^2)^3} \big)
 - \frac{1}{\sqrt{t_1}}(A_2|x|+D_2)(\f{1}{2\sqrt {t-t_1}}+\f{A_3 R^2}{\sqrt{t_1}(R^2-|x|^2)}+\theta)\f{2A_3R^2|x|}{\sqrt{t_1}(R^2-|x|^2)^2}\\
 & +2(\f{1}{2\sqrt {t-t_1}}+\f{A_3R^2}{\sqrt{t_1}(R^2-|x|^2}))^3 \geq\\
 &-\e \big( \f{2A_3R^2}{\sqrt{t_1}(R^2-|x|^2)^2}+\f{8A_3R^4}{\sqrt{t_1}(R^2-|x|^2)^3}\big)
 -\frac{1}{\sqrt{t_1}}(A_2R+D_2)(\f{1}{2\sqrt {t-t_1}}+\f{A_3 R^2}{\sqrt{t_1}(R^2-|x|^2)}+\theta)\f{2A_3R^3}{\sqrt{t_1}(R^2-|x|^2)^2}\\
& +( \f{3}{\sqrt {t-t_1}})\f{A_3^2R^4}{t_1(R^2-|x|^2)^2}+2\f{A_3^3R^6}{t_1\sqrt{t_1}(R^2-|x|^2)^3},
  \end{array}
$$
where we have used that $|x|\leq R$. One can verify that the r.h.s. of the above equality, for $R>1$, $\e\leq 1$, $t_1\leq 1$ and $A_3=A_3(T)$ large enough, is strictly positive. Therefore, $z$ is a strict supersolution of (\ref{eq:p}) in $(t_1,T]\times B_R(0)$ and for $\e\leq 1$.
\\

\bigskip
 
We next prove that 
$$
|p (t,x)|\leq z(t,x), \quad \text{in $(t_1,T]\times B_R(0)$}.
$$
To this end, we notice that $z(t,x)$ goes to $+\infty$ as $|x|\to R$ or as $t\to t_1$. Therefore, $|p|(t,x)-z(t,x)$ attains its maximum at an interior point of $(t_1,T]\times B_R(0)$. We choose $t_m\leq T$ such that the maximum of $|p|(t,x)-z(t,x)$ in the set $(t_1,t_m]\times B_R(0)$ is equal to $0$. If such $t_m$ does not exist, we are done. Let $x_m$ such that $ |p|(t,x)-z(t,x) \leq  |p|(t_m,x_m)-z(t_m,x_m)=0$ for all $(t,x) \in (t_1,t_m) \times B_R(0)$. At such point, we have
$$
0\leq \p_t (|p|(t_m,x_m)- z(t_m,x_m)),\quad  0\leq  -\Delta (|p|(t_m,x_m)- z(t_m,x_m)),
$$
$$
|p |(t_m,x_m)\nabla |p|(t_m,x_m)=z(t_m,x_m)\nabla z(t_m,x_m).
$$
Combining the above properties with the facts that $|p|$ and $z$ are respectively sub and strict supersolution of (\ref{eq:p}), we obtain that
$$
2(|p |(t_m,x_m)-\theta)^3-2(z(t_m,x_m)-\theta)^3<0.
$$
It follows that
$$
|p|(t_m,x_m) <z(t_m,x_m),
$$
which is in contradiction with the choice of $(t_m,x_m)$. We deduce that
$$
|p(t,x)|\leq z(t,x)=\f{1}{2\sqrt{t-t_1}}+\f{A_3(T)R^2}{\sqrt{t_1}(R^2-|x|^2)}+\theta(T), \quad \text{in $(t_1,T]\times B_R(0)$}.
$$
The above equality holds for all $R>1$. We let $R\to \infty$ to obtain
$$
|p (t,x)|\leq \f{1}{2\sqrt{t-t_1}}+ \frac{A_3(T)}{\sqrt{t_1}} +\theta(T),  \quad \text{in $(t_1,T]\times \R^d$}.
$$
It follows that 
$$
|p (t,x)|\leq   \frac{A_4(T)}{\sqrt{t_1}} +\theta(T),  \quad \text{in $(2t_1,T]\times \R^d$}.
$$
Finally, choosing $t_1=\frac{t_0}{2}$ we obtain (\ref{1}).
\bigskip
 
We conclude by providing the proof of  Lemma \ref{lem:lb}: 

\begin{proof}[Proof of Lemma  \ref{lem:lb}]
We first notice thanks to (\ref{eq.5b}) that $n_\e$ satisfies
 $$
-K_2n_\e \leq \e\p_t n_\e-\e^2\Delta n_\e.
 $$ 
Using the heat kernel and assumption (\ref{as:n0-below}), we obtain that
 $$
\f{\e^{\f d2}}{(4\pi t)^{\f d 2}}\int_{|y-x_0|\leq L_0} e^{-\f{(x-y)^2}{4\e t }-\f{M_0+K_2t}{\e}}  dy \leq n_\e(t,x), \quad \forall (t,x)\in \R^+\times \R^d.
$$
We deduce that
 $$
\f{\e^{\f d2}}{(4\pi t)^{\f d 2}} |B_{L_0}(x_0)|  e^{-\f{2|x|^2+2(L_0+|x_0|)^2}{4\e t }-\f{M_0+K_2t}{\e}}  \leq n_\e(t,x), \quad \forall (t,x)\in \R^+\times \R^d,
$$
and hence
$$
\e \log (\f{\e^{\f d2}}{(4\pi t)^{\f d 2}} |B_{L_0}(x_0)| ) - \f{|x|^2+(L_0+|x_0|)^2}{2 t } -(M_0+K_2t) \leq u_\e(t,x), \quad \forall (t,x)\in \R^+\times \R^d.
$$
It follows that
$$
\e \log (\f{\e^{\f d2}}{(4\pi t)^{\f d 2}} |B_{L_0}(x_0)| ) - \f{|x|^2+(L_0+|x_0|)^2}{2 t_1 } -(M_0+K_2t) \leq u_\e(t,x), \quad \forall (t,x)\in (t_1,+\infty)\times \R^d.
$$
Finally (\ref{low-b-u}) follows for $\e\leq 1$, choosing constants $A_1$, $B_1$ and $D_1$ large enough and noticing that $\log(t)$ goes more slowly that $t$ to the infinity.
\end{proof}
%--------------------------------------

\section{Lipschitz bounds for equations (\ref{eq.i1})-(\ref{eq.i2})}
\label{ap3}

%--------------------------------------

Here we prove that $u_\epsilon$ are locally uniformly Lipschitz without assuming that the latter are differentiable. The proof follows the same ideas as in section \ref{lipx}.\\

Let $\overline{c}= \frac{2 L_1 b_M}{b_m}$. From (\ref{eq.iH1}) we have

\begin{align}
\nonumber&\partial_t \big(u_\epsilon(t,x+h)-u_\epsilon(t,x)+\overline{c}h\big(2u_\epsilon(t,x+h)-u_\epsilon(t,x)\big)\big)-(1+2\overline{c}h)R(x+h,I_\epsilon)+(1+\overline{c}h)R(x,I_\epsilon)\\
\nonumber&=\int K(z)b(x+h+\epsilon z,I_\epsilon)e^{\frac{u_\epsilon(t,x+h+\epsilon z)-u_\epsilon(t,x+h)}{\epsilon}}dz
- \int K(z)b(x+\epsilon z,I_\epsilon)e^{\frac{u_\epsilon(t,x+\epsilon z)-u_\epsilon(t,x)}{\epsilon}}dz\\
\nonumber& +\overline{c}h\big(\int K(z)2b(x+h+\epsilon z,I_\epsilon)e^{\frac{u_\epsilon(t,x+h+\epsilon z)-u_\epsilon(t,x+h)}{\epsilon}}dz- \int K(z)b(x+\epsilon z,I_\epsilon)e^{\frac{u_\epsilon(t,x+\epsilon z)-u_\epsilon(t,x)}{\epsilon}}dz \big)
\end{align}

Define $\alpha=\frac{u_\epsilon(t,x+\epsilon z)-u_\epsilon(t,x)}{\epsilon}$, $\beta=\frac{u_\epsilon(t,x+h+\epsilon z)-u_\epsilon(t,x+h)}{\epsilon}$, $\varDelta(t,x)=2u_\epsilon(t,x+h)-u_\epsilon(t,x)$ and $w_\epsilon(t,x)=\frac{u_\epsilon(t,x+h)-u_\epsilon(t,x)}{h}+\overline{c}\varDelta(t,x)$. Using the convexity inequality 
$$e^\beta\leq e^\alpha+e^\beta(\beta-\alpha),$$
we deduce 
\begin{align}
\nonumber h&\partial_t w_\epsilon(t,x)-(1+2\overline{c}h)R(x+h,I_\epsilon)+(1+\overline{c}h)R(x,I_\epsilon)\\
\nonumber&\leq\int K(z) b(x+h+\epsilon z,I_\epsilon)\big(e^\alpha+e^\beta(\beta-\alpha) \big)dz- \int K(z)b(x+\epsilon z,I_\epsilon)e^\alpha dz\\
\nonumber& +\overline{c}h\big(\int 2K(z)b(x+h+\epsilon z,I_\epsilon)e^\beta dz- \int K(z)b(x+\epsilon z,I_\epsilon)e^\alpha dz \big)\\
\nonumber& \leq \int K(z)\big(b(x+h+\epsilon z,I_\epsilon)-b(x+\epsilon z,I_\epsilon)\big)e^\alpha dz\\
\nonumber&+ \int K(z) b(x+h+\epsilon z,I_\epsilon)e^\beta\big(\beta-\alpha+ \overline{c}h \frac{\varDelta (t,x+\epsilon z)-\varDelta(t,x)}{\epsilon}\big)dz\\
\nonumber&
+\overline{c}h\int K(z)b(x+h+\epsilon z,I_\epsilon)e^\beta(2-2\beta+\alpha)dz-\overline{c}h\int K(z)b(x+\epsilon z,I_\epsilon)e^\alpha dz.
\end{align}

From assumptions (\ref{eq.5b}) and (\ref{eq.6i}) it follows that 

\begin{align}
\nonumber\partial_t w_\epsilon(t,x)&\leq \int K(z) b(x+h+\epsilon z,I_\epsilon)e^\beta\frac{w_\epsilon(t,x+\epsilon z)-w_\epsilon(t,x)}{\epsilon}dz\\
\nonumber& +K_2+3\overline{c}K_2+\int K(z)\big( \overline{c}b_Me^\beta(2-2\beta+\alpha)+(L_1b_M-\overline{c}b_m)e^\alpha\big) dz.
\end{align}

Notice that 
$$ \overline{c}b_Me^\beta(2-2\beta+\alpha)+(L_1b_M-\overline{c}b_m)e^\alpha= \overline{c}b_Me^\beta(2-2\beta+\alpha)-L_1b_M e^\alpha,$$
is bounded from above. Indeed if we first maximize the latter with respect to $\beta$ and then with respect to $\alpha$ we obtain

$$\overline{c}b_Me^\beta(2-2\beta+\alpha)-L_1b_M e^\alpha\leq 2\overline{c}b_Me^{\frac{\alpha}{2}}-L_1b_M e^\alpha\leq \frac{b_M \overline{c}^2 }{L_1}.$$

We deduce 
\begin{align}
\nonumber\partial_t w_\epsilon(t,x) &\leq \int K(z) b(x+h+\epsilon z,I_\epsilon)e^\beta\frac{w_\epsilon(t,x+\epsilon z)-w_\epsilon(t,x)}{\epsilon}dz +G,
\end{align}
where $G$ is a constant. Therefore by the maximum principle, (\ref{eq.iH3}) and (\ref{lowerb}), we have

$$w_\epsilon(t,x)\leq Gt+\parallel\nabla u_\epsilon^0\parallel_{L^\infty}-2\overline{c}A|x+h|+2\overline{c}B-\overline{c}u_\epsilon^0(x=0)+\overline{c}\parallel \nabla u_\epsilon^0\parallel_{L^\infty}|x|.$$

Using again (\ref{eq.iH3}) and (\ref{lowerb}) we conclude that

\begin{align}
\frac{u_\epsilon(t,x+h)-u_\epsilon(t,x)}{h}&\leq (G+2\overline{c}K_2)t+\overline{c}\big(-A+\parallel\nabla u_\epsilon^0\parallel_{L^\infty}\big)\big(|x|+2|x+h|\big)\\
\nonumber&+3\overline{c}B+\parallel\nabla u_\epsilon^0\parallel_{L^\infty}-3\overline{c}\inf u_\epsilon^0(x=0).
\end{align}

\end{document}